\documentclass[12point]{amsart}
\usepackage{geometry}                                   \usepackage[parfill]{parskip}     
\usepackage{amsthm}
\newtheorem*{lemma}{Lemma}
\newtheorem*{CT}{Cobham's Theorem}

\usepackage{amsmath}

\def \sds  {symbolic dynamical system}
\def \sms {substitution minimal system}

\def \sub  {substitution}
\def \topconj {topologically conjugate}

\def \st   {such that }

\def  \s {\sigma}

\def \ZZ {\mathbb Z}
\def \th {\theta}

\usepackage{graphicx}
\usepackage{amssymb}
\usepackage{epstopdf}
\title{A short proof of a theorem of Cobham on substitutions}

\author{Ethan M. Coven}
\address{Ethan M. Coven,
Department of Mathematics,
Wesleyan University,
Middletown, CT 06459}
\email{ecoven@wesleyan.edu}
\author{Andrew Dykstra}
\address{Andrew Dykstra,
Department of Mathematics,
Hamilton College,
Clinton, NY 13323}
\email{adykstra@hamilton.edu}
\author{Michelle LeMasurier}
\address{Michelle LeMasurier,
Department of Mathematics,
Hamilton College,
Clinton, NY 13323}
\email{mlemasur@hamilton.edu}

\date{August 17, 2011}

\begin{document}

\LARGE

\maketitle

\begin{abstract}
This paper is concerned with the lengths of constant length \sub s that generate \topconj\ systems.  We show that if the systems are infinite, then these lengths must be powers of the same integer.  This result is a dynamical formulation of a special case of a 1969 theoretical computer science result of Alan Cobham~\cite{Cob}.  Our proof is rather simple.
\end{abstract}

\section{Introduction}

This paper is concerned with the lengths of constant length \sub s that generate \topconj\ systems.  We show that if the systems are infinite, then these lengths must be powers of the same integer.  This result is a dynamical formulation of a special case 
(for constant length \sub s) of a 1969 theoretical computer science result of Alan~Cobham \cite{Cob}.  Our proof is rather simple.

Fabien Durand (Theorem ~8 of \cite{Dur}) put Cobham's Theorem in the setting of symbolic dynamics.
For a brief discussion of the history of Cobham's Theorem and how it came to symbolic dynamics via computer science and logic, see \cite{Dur}.  Included there are S. Eilenberg's remarks on the "highly technical" nature of Cobham's original proof \cite{Eil}, as well as comments on subsequent proofs.

Durand's proof, unlike ours, works for non-constant length \sub s, too, and hence is more complicated.  
Our proof is, we believe, the simplest proof yet of (this special case of) Cobham's Theorem.  It relies on the characterization in \cite{CKL} of symbolic minimal systems topologically conjugate to constant length \sms s.

The authors thank Fabien Durand for several enlightening electronic discussions about substitutions and Cobham's Theorem.

\section{Basic Concepts}

A {\it dynamical system\/} is a pair $(X,T)$, where
$T : X \to X$ is a homeomorphism.  Dynamical systems
$(X,T)$ and $(Y,S)$ are considered ``the same'' if the actions of~$T$ on~$X$ and of~$S$ on~$Y$ are the same, only the names of the points have been changed.
Formally, $(X,T)$ and $(Y,S)$ are {\it \topconj\/}
iff there is a homeomorphism $\varphi : X \to Y$,
called a {\it topological conjugacy\/},
\st $\varphi \circ T \equiv S \circ \varphi$.
A dynamical system $(X,T)$ is called {\it minimal\/}
iff $X$ contains no nonempty, closed, $T$-invariant subset.

A {\it symbolic dynamical system\/} is a
dynamical system $(X,\s)$, where $\s$ is the (left) shift and $X$ is a closed, shift-invariant subset of some
$A^\ZZ = \prod_{-\infty}^{\infty} A$,
the space of all doubly infinite sequences 
with entries from the finite {\it alphabet\/}
~$A$.
Here $A$ has the discrete topology and $A^\ZZ$
the product topology.  (If $\# A\ge 2$, then $A^\ZZ$ is homeomorphic to the Cantor set.)
The {\it shift\/} $\s : A^\ZZ \to A^\ZZ$ is defined by $[\s(x)]_i := x_{i+1}$ for every $i \in \ZZ$.
To avoid notational clutter, when the domain is clear we will use the same symbol~$\s$ to denote the shift on all~$A^\ZZ$ and on all closed, shift-invariant subsets.

If $X$ is a closed, $\s^k$-invariant subset of some 
$A^\ZZ$, then $(X,\s^k)$ is
\topconj\ to $(Y,\s)$, where the symbols of~$Y$
are the words of length~$k$ of~$X$ that appear starting at places that are multiples of~$k$.  Thus $(X,\s^k)$ is
a symbolic dynamical system.

A {\it \sub\/} of constant length $L \ge 2$ is a mapping
$\th : A \to A^L$, the words of length~$L$, where $A$
is a finite alphabet.  The most famous example is the Morse (Morse-Thue-Prouhet) substitution  
$0 \mapsto 01, 1 \mapsto 10$.

A \sub\ $\th : A \to A^L$ maps $A^2$ to $A^{2L}$ by juxtaposition: $\th(ab):= \th(a)\th(b).$
In the same way, for every $k \ge 2$, 
$\th$ maps~$A^k$ to~$A^{kL}$,
and $A^{\ZZ}$ to itself.
The \sub\ $\th$ is {\it primitive\/}  
iff there exists~$i$ \st
$A \subseteq \th^i(a)$ for every $a \in A$.
If  $\th$ is primitive,
then there is a unique smallest \sds\ $(X_\th,\s)$
\st every word $\th^i(a)$, $a \in A$, $i \ge 1$
appears in~$X_\th$.
Then $(X_\th,\s)$ is minimal and is
called the {\it \sms\ generated by}~$\th$.

Finally, it is clear that for every $n \ge 2$, $\th$ and $\th^n$ generate the same \sms.

\section{Cobham's Theorem}

\begin{CT}
The lengths of primitive, constant length substitutions that generate topologically conjugate infinite \sms s are powers of the same integer.
\end{CT}

\proof Suppose that $\theta$ and $\zeta$ are primitive, constant length substitutions that generate topologically conjugate infinite \sms s $(X_{\theta},\sigma)$ and $(X_{\zeta},\sigma)$.
In \cite{Dek} F.~M.~Dekking proved that    
the set of prime divisors of the length of a primitive, constant length substitution that generates an infinite \sms\ is  a topological conjugacy invariant of that  \sms.  We show that if the lengths of ~$\th$ and~$\zeta$ are not powers of the same integer, then both
$(X_\th,\s)$ and $(X_\zeta,\s)$ are \topconj\ to a \sms\ generated by a \sub\ whose length has fewer prime factors than do the lengths of ~$\th$ and~$\zeta$, contradicting Dekking's Theorem.

So let the lengths of $\theta$ and $\zeta$ be

$$ p_1^{m_1}p_2^{m_2}\cdots p_k^{m_k}
\text{   and   }
p_1^{n_1}p_2^{n_2}\cdots p_k^{n_k},
$$

\noindent where the $p$'s are distinct primes and the $m$'s and  $n$'s are positive.

Let $J$  be such that
\begin{equation}
\frac{m_J}{n_J} \le \frac{m_i}{n_i}
\text{ for all }i.
\end{equation}  
If these lengths are not powers of the same integer, then
\begin{equation}\frac{m_J}{n_J} < \frac{m_i}{n_i}
\text{ for some } i.
\end{equation}

Replace  $\theta$ by $\theta^{n_J}$ and 
$\zeta$ by $\zeta^{m_J}$. The powers generate the same \sms s as do $\theta$ and $\zeta$.  
It follows from~(1) and~(2)
that the lengths $M$ and $N$ of
$\theta^{n_J}$ and $\zeta^{m_J}$
satisfy $M=RN$, where $R > 1$ and
has fewer prime factors than do $M$ and ~$N$.
By the Lemma  below, there is a \sms\ generated by a primitive substitution of constant length~$R$ that is topologically conjugate to
$(X_{\theta},\sigma)$ and $(X_{\zeta},\sigma)$.  This contradicts Dekking's Theorem.
\qed

\begin{lemma}
Suppose that $\theta$ and $\zeta$ are primitive  \sub s of constant lengths $M$ and~$N$ that generate topologically conjugate infinite \sms s $(X_{\theta},\sigma)$ and $(X_{\zeta},\sigma)$.  If $M=RN$, where $R>1$, then there is a primitve substitution of constant length~$R$ that generates a \sms\ that is topologically conjugate to $(X_{\theta},\sigma)$ and $(X_{\zeta},\sigma)$.
\end{lemma}

\begin{proof}
Replace $X_{\zeta}$ by its image in ~$X_\th$ under the topological conjugacy from 
$(X_{\zeta},\sigma)$ to $(X_{\theta},\sigma)$.
Since the condition of the theorem on p.~1450 of~\cite{CKL} is preserved by topological conjugacy,
there are proper $\s^M$- and $\s^N$-invariant subsets $X_M$ and $X_N$ of ~$X_\th$ \st $(X_M,\s^M)$ and
$(X_N,\s^N)$ both are \topconj\ to $(X_\th,\s)$.

Since
$$\bigcup_{i=0}^{N-1} \s^i
[X_M \cup \s^N(X_M) \cup \cdots \cup \s^{(R-1)N}(X_M)] = X_\th,
$$

it follows from the minimality of~$(X_N,\s^N)$
that one of the members of the (large) union is
~$X_N$.  Thus $X_N$ has a proper subset,
$\s^i(X_M)$ for some~$i$, 
\st $(X_N,\s^N)$ is topologically conjugate to
$(\s^i(X_M),(\s^N)^R)$.

Therefore, by the theorem in ~\cite{CKL},
$(X_N,\s^N)$ is \topconj\ to a \sms\ generated by a primitive \sub\ of constant length~$R$.  Hence so is 
$(X_\th,\s)$.
\end{proof}

\end{document}